\documentclass[10pt]{article}
\usepackage{amsmath, amssymb, amsthm, pb-diagram, euscript, mathrsfs}

\usepackage[T1]{fontenc}
\usepackage[sc]{mathpazo}
\usepackage[usenames]{color}
\usepackage{colortbl}

\DeclareFontFamily{T1}{pzc}{}
\DeclareFontShape{T1}{pzc}{m}{it}{1.8 <-> pzcmi8t}{}
\DeclareMathAlphabet{\mathpzc}{T1}{pzc}{m}{it}

\title{Noncommutative Fr{\'e}chet Spaces and Unbounded Bivariant K-Theory }
\author{Nikolay P. Ivankov\\
Max-Planck-Institut f\"ur Mathematik, Bonn}

\theoremstyle{plain}
\newtheorem{prop}{Proposition}[section]

\newtheorem{lem}[prop]{Lemma}
\newtheorem{cor}[prop]{Corollary}
\newtheorem{thm}[prop]{Theorem}

\theoremstyle{definition}
\newtheorem{defn}[prop]{Definition}
\newtheorem{exm}[prop]{Example}
\newtheorem{rem}[prop]{Remark}

\theoremstyle{remark}

\newcommand{\A}{\mathcal{A}}

\newcommand{\Hc}{\mathcal{H}}

\newcommand{\Lc}{\mathcal{L}}

\newcommand{\gtb}{\mathfrak{b}}

\newcommand{\Cb}{\mathbb{C}}
\newcommand{\Kb}{\mathbb{K}}
\newcommand{\Nb}{\mathbb{N}}

\newcommand{\Zb}{\mathbb{Z}}

\newcommand{\Asm}{\mathscr{A}}
\newcommand{\Bsm}{\mathscr{B}}

\newcommand{\Lsm}{\mathscr{L}}

\newcommand{\Te}{\Theta}

\newcommand{\Om}{\Omega}

\newcommand{\al}{\alpha}
\newcommand{\be}{\beta}
\newcommand{\gm}{\gamma}

\newcommand{\la}{\lambda}

\newcommand{\sg}{\sigma}

\newcommand{\diag}{\mathrm{diag}}

\newcommand{\End}{\mathrm{End}}

\newcommand{\ad}{\mathrm{ad}}

\newcommand{\Qin}[1]{\Psi^{({#1})}}
\newcommand{\KK}{\mathrm{KK}}
\newcommand{\Ai}[1]{\A^{({#1})}}

\newcommand{\iy}{\infty}

\newcommand{\Dn}{{\mathop{D\!\!\!\!/}}}

\newcommand{\ord}{\mathrm{ord}}

\newcommand{\nix}{\phantom{.}}
\newcommand{\btr}[1]{{#1}(1+{#1}^2)^{-\frac{1}{2}}}

\begin{document}

\maketitle
\begin{abstract} In this paper we introduce an abstract approach to higher order smooth systems on $C^*$-algebras in contest of Baaj-Julg picture of $\KK$-theory.
\end{abstract}

\section{Introduction}
The unbounded picture of $\KK$-theory was proposed by Saad Baaj and Pierre Julg in their paper \cite{BaJu}, published in 1983. We recall the definition of unbounded $\KK$-cycle.
\begin{defn}
Let $A$, $B$ be $C^*$-algebras. An \emph{unbounded $\KK$-cycle} over $(A,B)$ is a set of data $(E,\pi,D)$, where $E$ is a $C^*$-bimodule, $\pi$ a $C^*$-algebra representation $\pi\colon A\to\End_B(E)$ and $D$ an unbounded selfadjoint regular $B$-linear operator on $E$ \cite{Woro}, satisfying 
\begin{enumerate}
\item
Operators $\pi(a)D(1+D^2)^{-\frac{1}{2}}\in\Kb_B(E)$ for all $a\in A$
\item 
The set of $a\in A$ for which $[D,\pi(a)]$ extends to a bounded operator on $E$ is dense in $A$. 
\end{enumerate}
\end{defn}

Baaj and Julg have proved that this approach significantly simplifies the calculations of the representative $\KK$-cycle in the case when one considers the outer product in $\KK$-theory.

However, the most challenging question of the applicability of Baaj-Julg approach to the Kasparov product still remains open. A sufficient advance in this direction was made Dan Kucerovsky in \cite{Kuce} and later by Bram Mesland in \cite{Mesl}. The goal of the latter work is to show that the calculation of the representative of Kasparov product of two $\KK$-cycles may be reduced from the Kasparov's technical result to a simple formula. The burden of calculations involved in the technical result will then be reloaded on the conditions imposed on the dense subspaces of $C^*$-algebras and $C^*$-modules and unbounded operators involved in the calculations. These properties, in turn, are believed to be easily checked for interesting spaces such as real spectral triples \cite{ConNCG},\cite{Vari}.

One of the main differences between classical Baaj-Julg approach and the one by Mesland is that one considers an analogue of $C^1$ functions in the first one, whenever for the construction of Kasparov product in the way described in \cite{Mesl} one needs at least $C^2$-algebras generated by $D$ to be dense in the original $C^*$-algebra. In addition, in the case of the algebras of functions on the manifolds, we have natural $C^\iy$-structures which have some useful properties (such as nuclearity of $C^\iy(M)$). Therefore the study of $C^n$-subalgebras of $C^*$-algebras may become an important branch of development is the unbounded $\KK$-theory. 

\section{Definitions}

Unlike the commutative setting, where one has a canonical way to construct the subalgebras of $C^n$ functions my means of the derivations on local atlases, noncommutative setting does not a priori have such a structure. There have already been proposed several approaches to this question. Among the most remarkable ones are the approach of noncommtative geometry and the one by Blackadar and Cuntz. In the first one the algebra of smooth elements for the $C^*$-algebra is actually defined via the additional structure of spectral triple. The Blackadar-Cuntz approach is, in contrast, very abstract and uses only the Banach algebra structures, having nothing to do with any unbounded operators.

In what following, we are going to combine these two.

\begin{defn}
Let $A$ be a $C^*$-algebra. An \emph{$n$-smooth system} (or $C^n$-system) $\Asm$ on $A$ is an inverse system of pre-$C^*$-subalgebras of $A$
$$
\Ai{n}\subseteq \Ai{n-1} \subseteq \dots \subseteq \Ai{1} \subseteq \Ai{0} = A
$$
such that all $\Ai{j}$ are also isomorphic to operator algebras, which we shall abusively denote by $\Ai{j}$, and the operator algebra maps induced by the inclusions $\Ai{k}\to\Ai{k-1}$ are completely bounded for $k=1,\dots,n$. It will be called \emph{$\iy$-smooth} (or $C^\iy$) if the system of these subalgebras is infinite and the inverse limit of this system $\Ai{\iy}$ is also a pre-$C^*$-algebra. The number $n$ (including the case $n=\iy$) will be called the \emph{order} of the smooth system and will be denoted as $\ord(\Asm)$. We will also employ the notation $\|\cdot\|_n$ for the operator norm on $\Ai{n}$.
\end{defn}

Though in this definition we use operator algebras instead of Banach algebras, this approach is in some sense even more general then the one in \cite{BlaCu}. We are now going to develop the framework that will relate it to the unbounded $\KK$-theory.

\begin{defn}
For given two $C^*$-algebras $A$ and $B$, an unbounded $(A,B)$-$\KK$-cycle $(E,D)$ and a natural number $k$ the \emph{fr{\'e}cetization} is a map $\mu \colon (A,B,E,D,k) \to \Ai{k}_{\mu,D}$, where $\Ai{k}_{\mu,D}$ is an operator algebra isomorphic (in pure algebraic sense) to a subalgebra of $A$.
\end{defn}

\begin{defn}
A \emph{smooth system on an algebra $A$ generated by operator $D$ with respect to the fr{\'e}chetization $\mu$} is the longest sequence of nested subalgebras $\Ai{k}_{\mu,D}$ of $A$ with the starting point $\Ai{0}=A$ satisfying the conditions of smooth system. We will denote this system by $\Asm_{\mu, D}$
\end{defn}

\begin{defn}
Let $\Asm$ be a smooth system on a $C^*$-algebra $A$ with $\ord{\Asm}\ge n$, $n\in\Nb$, and let $B$ be another $C^*$-algebra. We shall say that the unbounded $(A,B)$-$\KK$-cycle $(E,D)$ is \emph{$n$-smooth ($C^n$) with respect to $\Asm$} if $\Ai{k}\subset \Ai{k}_{\mu, D}$ for all $k\le n$, and the inclusion morphism induces a completely bounded homomorphism of operator (pseudo)algebras. The set of all such cycles will be denoted by $\Qin{n}_\mu(\Asm,B)$. We say that $(E,D)\in\Qin{n}_{\mu}(\Asm,B)$ if $\ord(\Asm)=\iy$ and the same property holds for all $n\in\Nb$. Note that in this case $\Ai{\iy}_{\mu,D}$ will automatically be a pre-$C^*$-algebra. We shall use the notation $\|\cdot\|_{n,\mu,D}$ for the operator norm on $\Ai{n}_{\mu,D}$.
\end{defn}

Analogously to the definition of usual $\KK$ cycles we may define the set $\Qin{n}_{0,\mu}(\Asm,B)$ ($\Qin{n}_{1,\mu}(\Asm,B)$) for the cases when the cycle $(E,D)$ is even (resp. odd). 

The completely bounded maps between the $C^n$-algebras play an important of preservation of an important class of so-called smooth modules. We refer to \cite{Mesl} for an explicit picture. 

\begin{rem}
In fact the $\KK$-cycles, both bounded and unbounded, are defined as triples $(E,\pi,F)$ (resp. $(E,\pi,D)$), where $\pi\colon A \to \End^*_B(E)$ is a representation of $A$. However, to save the space, we are going to use the same symbol $a$ for the element of (a subalgebra of) the $C^*$-algebra $A$, an element of an operator algebra, and also $\pi(a)$. The author believes that such an identification would not leave to misunderstandings.
\end{rem}

\section{Relation to Classical $\KK$-Theory}

\begin{defn}
Let $\mu$ be a fr{\'e}chetization. We shall call $\mu$ \emph{commutator bounded} if for all $n\in \Nb$ and there exist positive numbers $C_n$ such that for any $C^*$-algebra $A$ and for any $a \in \Ai{n}_{\mu,D}$ one has that 
$$
\|a\|_{n,\mu,D} \le C_n\max\{\|a\|, \|\ad_D(a)\|, \ad^2_D(a)\|, \dots, \|\ad^n_D(a)\|\}
$$
Here $\ad^n_D(a) = [D;[\dots[D;a]\dots]]$, the $n$'th commutator of $D$ with $a$.
We include the case when  $\ad^k_D(\pi(a))$ does not extend to a bounded operator on $E$, setting $\|\ad^k_D(\pi(a))\|=\iy$.

The fr{\'e}chetization will be called \emph{analytic} if for any smooth system $\Asm_{\mu,D}$ the norm $\|\cdot\|_{n,\mu,D}$ is analytic with respect to $\|\cdot\|_{n-1,\mu,D}$ for all $1\le n \le \ord(\Asm_{\mu,D})$.
\end{defn}

We recall (cf. \cite{BlaCu},\cite{Mesl}), that if $\A$ is an algebra with the Banach norm $\|\cdot\|_\al$, and $A_\al$ is its closure with respect to this norm, then a Banach norm $\|\cdot\|_\be$ on $\A$ is called \emph{analytic} with respect to $\|\cdot\|_\al$, if for all $a\in\A$ such that $\|a\|_\al<1$ holds
$$
\limsup_{n\to\iy}\frac{\ln\|a^n\|_\be}{n}\le 0
$$

The main consequence of analyticity of one norm with respect to another is the stability of $A_\be$ with respect to the holomorphic functional calculus on $A_\al$ (cf. \cite{BlaCu},\cite{Mesl}). Here $A_\be$ is the completion of $\A$ with respect to $\|\cdot\|_\be$.

Observe also, that if $\|\cdot\|_\gm \le C\|\cdot\|_\be$ then it is also analytic with respect to $\|\cdot\|_\al$. Indeed,
$$
\limsup_{n\to\iy}\frac{\ln\|a^n\|_\gm}{n}\le\limsup_{n\to\iy}\frac{\ln C\|a^n\|_\be}{n}= \limsup_{n\to\iy}\frac{\ln C+\ln\|a^n\|_\be}{n} \le 0
$$

Now we are ready to formulate the main result of the section.

\begin{thm}\label{smth_alg_thm_relation_to_KK}
Let $\mu$ be a fr{\'e}chetization, which is commutator bounded and analytic. Then, for any separable $C^*$-algebra $A$ and any set of isomorphism classes of $C^*$-alberas $\Lambda$ there is such an $\iy$-smooth system $\Asm$ on $A$ that for any $B\in\Omega$ there is a surjective map $\Qin{\iy}(\Asm, B) \to \KK(A,B)$, induced by the bounded transform map.
\end{thm}

Before we proceed to the proof, we shall discuss some more specific formulations, allowing us to apply the result in more concrete situations. First of all it should be noted that the notion of the fr{\'e}chetization was introduced by the author because he has encountered different ways to define the $C^n$ algebra by means of an unbounded $\KK$-cycle. Two of them, which are now playing the most important role in the theory, are discussed in the next section, and both satisfy the conditions of the theorem with $C_n = 2^n$.  

As to the set $\Lambda$ mentioned in the formulation, we may have the following examples.

\begin{exm}
Let $\Lambda=\{\Cb\}$. Then the unbounded $\KK$-cycles in $\Qin{n}(\Asm,B)$ are spectral triples in the widest sense. It should be noted that the sense is indeed wide: as we shall see later, the smooth system $\Asm$ may be very far from the familiar one.
\end{exm}

\begin{exm}
It is a well-known fact that every $C^*$-algebra may be represented as a $C^*$-subalgebra of $\Bsm(\Hc)$ for a separable Hilbert space $\Hc$. Therefore the isomorphism classes of separable $C^*$-algebras form a set. Thus for any separable $C^*$-algebra $A$ we may choose a unique smooth system $\Asm$, such that for any separable $C^*$-algebra $B$ there will be a surjective map $\Qin{n}(\Asm,B)\to\KK(A,B)$. 
\end{exm}

In order proceed to the proof of the theorem \ref{smth_alg_thm_relation_to_KK} we first need to prove following lemmas.

\begin{lem}\label{smth_alg_lem_tech_lemma}
Let $A$ be a separable $C^*$-algebra. Then for any $C^*$-algebra $B$ and any element $[(E,F)]\in\KK(A,B)$ there exists an unbounded $(A,B)$-$\KK$-cycle $(E,D)$, such that $[E,\btr{D}]=[(E,F)]$ and the set of such $a\in A$, that $\ad_D^n(a)$ extends to a bounded operator on $E$, is dense in $A$.
\begin{proof}
This result is a generalization of the Theorem 17.11.4 form \cite{Bla}. Fix a total system $\{a_j\}$ of $A$. For given $F$ there exists a strictly positive element $h\in \Kb(E)$ of degree $0$ which commutes with $F$ \cite{Bla}. Now, according to \cite[3.12.14]{Ped1} there exists an approximate unit $u_k$ for $\Kb(E)$, contained in $C^*(h)$, quasicentral for $A$, with the property that $u_k\ge 0$, $u_{k+1}\ge u_k$ and $u_{k+1}u_{k}=u_{k}$ for all $k\in\Nb$. Denote $d_k=u_{k+1}-u_k$. Passing, when needed, to a subsequence, we may assume that $\|d_k[F;a_j]\|<2^{-k^2}$ and $\|[d_k;a_j]\|<2^{-k^2}$ for all $k\ge j+1$. Set $X=\widehat{C^*(h)}\approx \sg(h)\setminus\{0\}$ and let $X_n$ be the support of $u_k$. Then $\langle X_k\rangle$ is an increasing sequence of compact subsets of $X$ and $X=\bigcup_{k=1}^\iy X_k$. Put 
\begin{equation*}\label{smth_alg_eqn_r_construction}
r_k=\sum_{l=1}^k 2^l d_l
\end{equation*}

This sequence converges pointwise on $X$ to an unbounded function $r$. Observe
that $r\ge 2^{k}$ on $X\setminus X_k$, so that $R=r^{-1}$ defines an element of
$C^*(h)$. Note also, that $d_k$ defines a bounded function on the space X and, since $\|d_k\|\le 1$ and $d_k d_{k-l}=0$ for all $k\ge 3$ and $2\le l 
\le k-1$, we obtain that 
$$
\|r_k\|\le\max_{l=2,\dots,k}\{\|2^{l-1} d_{l-1} + 2^l d_l\|\} \le 3\cdot 2^{k-1} < 2^{k+1}
$$ 

Let now $D=Fr$. Then $D=D^*$ and $(1+D^2)^{-1}$ extends to
$R^2(1+R^2)^{-1}\in\Kb_A(E)$. 

Observe that, since $F$, $r_k$ and $d_k$ commute for all $k$,

\begin{equation*}
\begin{split}
&\ad_{Fr_{k+1}}^n(a_j)-\ad_{Fr_k}^n(a_j)=\\
&\ad_{F(r_k+2^{k+1}d_{k+1})}^n(a_j)-\ad_{Fr_k}^n(a_j)=\\
&\sum_{l=0}^n C^l_n\ad_{ Fr_k}^{n-l}\left(\ad_{F\cdot 2^{k+1}d_{k+1}}^l(a_j)\right)-\ad_{Fr_k}^n(a_j)=\\
&\sum_{l=1}^n 2^{l(k+1)} C^l_n\ad_{Fr_k}^{n-l} \left(\ad _{Fd_{k+1}}^l(a_j)\right)=\\
&\sum_{l=1}^n 2^{l(k+1)}C^l_n \ad_{ Fr_k}^{n-l} \left(\ad_{Fd_{k+1}}^{l-1}([Fd_{k+1};a_j])\right)=\\
&\sum_{l=1}^n 2^{l(k+1)}C^l_n \ad_{Fr_k}^{n-l} \left(\ad_{Fd_{k+1}}^{l-1} (F[d_{k+1};a_j]+[F;a_j]d_{k+1})\right)\\
\end{split}
\end{equation*}
where $C^l_n$ are binomial coefficients. Now since $\|F\|=1$, $\|d_{k+1}\|\le 1$
and $\|r_k\|< 2^{k+1}$, we obtain that $\|\ad_{ Fd_{k+1}}(b)\|\le 2\|b\|$
and $\|\ad_{Fr_k}(b)\|\le 2^{k+2}\|b\|$ for any bounded operator $b$, we estimate for $k\ge j+1$:
\begin{equation*}
\begin{split}
&\left\|\sum_{l=1}^n 2^{l(k+1)}C^l_n \ad_{Fr_k}^{n-l} \left(\ad_{F d_{k+1}}^{l-1} (F[d_{k+1};a_j]+[F;a_j]d_{k+1})\right) \right\|\le\\
&\sum_{l=1}^n 2^{l(k+1)}C^l_n\left\|\ad_{Fr_k}^{n-l} \left(\ad_{F d_{k+1}}^{l-1} (F[d_{k+1};a_j]+[F;a_j]d_{k+1})\right)\right\|
<\\
&\sum_{l=1}^n 2^{l(k+1)}C^l_n\cdot 2^{(k+2)(n-l)}\left\|\ad_{Fd_{k+1}}^{l-1} (F[d_{k+1};a_j]+[F;a_j]d_{k+1})\right\|\le\\
&\sum_{l=1}^n 2^{l(k+1)}C^l_n\cdot 2^{(k+2)(n-l)}\cdot 2^{l-1}\left\| F[d_{k+1};a_j]+[F;a_j]d_{k+1}\right\|\le\\
\end{split}
\end{equation*}
\begin{equation*}
\begin{split}
&\sum_{l=1}^n 2^{l(k+1)}C^l_n\cdot 2^{(k+2)(n-l)}\cdot 2^{l-1} \cdot (2^{-k^2}+2^{-k^2})=\\
&\sum_{l=1}^n C^l_n 2^{l(k+1)+(k+2)(n-l)+(l-1)+1-k^2}=\\
&\sum_{l=1}^n C^l_n 2^{kn+2n-k^2}=\\
&2^n \cdot 2^{kn+2n-k^2}=\\
&2^{-k^2+(k+3)n}
\end{split}
\end{equation*}

Thus we have that the sequence $\{\ad_{Fr_k}^n(a_j)\}$ is norm convergent. Hence the operator
$$
\ad_D^n(a_j)=\ad_{Fr}^n(a_j)
$$
extends to a bounded operator on $E$. Thus the set of all $a\in\A$ such that $\ad_D^n D(a)$ extends to a bounded operator on $E$ is dense in $\A$. Pointing out, that it is true for $n=1$ and observing that 
$$
D(1+D^2)^{-1/2}=F(1+R^2)^{-1/2}
$$ 
is a "compact perturbation" of $F$, we obtain that $(E,D)$ is an
unbounded $(A,B)$-$KK$-cycle and that $[(E,F)]=[(E,\gtb(D))]$. QED.
\end{proof}
\end{lem}

We have actually shown more then we have claimed in the formulation of Lemma \ref{smth_alg_lem_tech_lemma}. Namely, we proved that for any element $(E,F)$ and any total system $\{a_j\}$ one we may construct such $D$ that 
\begin{equation}\label{smth_alg_eqn_commutator_estimate}  
\|\ad^n_D(a_j)\|\le c_{n,j}
\end{equation}
where $C_{n,j}$ is a positive number that does not depend neither on the choice of $F$ nor on $\{a_j\}$. This observation lets us prove the next lemma.

\begin{lem}\label{smth_alg_lem_general_lemma_1}
Let $\mu$ be some fr{\'e}chetization, $A$ be a separable $C^*$-algebra, $\{a_j\}$ - an arbitrary total system on $A$ and $\Omega$ be a set of such unbounded $(A,B_\omega)$-$\KK$-cycles $(E_\omega,D_\omega)$ that
\begin{itemize}
\item
$\ord(\Asm_{\mu, D_\omega})=\iy$
\item
For each $n$ and the operator algebra $\Ai{n}_{n,\mu,D_\omega}$ we have that $\|a_j\|_{n,\mu, D_\omega} \le K_{n,j}$, where $K_{n,j}$ are some positive numbers independent of $(E_\omega,D_\omega)\in\Om$. 
\end{itemize} 
Then there is a nested system of dense subalgebras $\Ai{n}$ in $A$, satisfying all the properties of smooth system except, possibly, for holomorphic stability, with a property that $\Ai{n}\subseteq \Ai{n}_{\mu,D}$ and the map induced by this inclusion is completely bounded.
\begin{proof}

We iteratively define the matrix norms 
\[
\|(a)_{ik}\|_n := \max\{\sup_{\omega\in \Om} \|(a)_{ik}\|_{n,\mu,D};\,\nix_m \|(a)_{ik}\|_{n-1}\}
\]
with $\|a\|_0$ being the $C^*$-norm on $A$ and let $\Ai{n}$ be the completion of $\mathrm{span}(\{a_j\})$ with respect to $\|a\|_n$. We have that $\Ai{n}$ is dense in $A$, and also is an algebra. The above matrix norms are finite for all $(a)_{ik}\in M_n(\Ai{n})$, and one may check directly that they make $\Ai{n}$ into an $\Lsm^{\iy}$ matricially normed space (see, for instance, \cite{EfRu0}, \cite{EfRu} for definition). It is also easy to check that the multiplication on $\Ai{n}$ is completely contractive. Indeed, for $\Ai{0}=A$, so the claim holds for $n=0$. Suppose that it is true for $n-1$. Then for $n$ we have
\[
\begin{split}
&\|(a)_{ik}(b)_{pq}\|_n \le \\
& \le  \max\{\|(a)_{ik}(b)_{pq}\|_{n-1};\, \sup_{\omega\in\Om}  \|(a)_{ik}(b)_{pq}\|_{n,\mu,D_\omega}\} \\
&\le \max\{\|(a)_{ik}\|_{n-1}\nix_m\|(b)_{pq}\|_{n-1};\, \sup_{\omega\in\Om} \|(a)_{ik}\|_{n,\mu, D_\omega} \sup_{\omega\in\Om}  \|(b)_{pq}\|_{n,\mu,D_\omega}\} \\
& \le \max\{\|(a)_{ik}\|_{n-1};\, \sup_{\omega\in\Om}  \|(a)_{ik}\|_{n,\mu, D_\omega}\} \max\{\|(b)_{pq}\|_{n-1}; \, \sup_{\omega\in\Om}  \|(b)_{pq}\|_{n,\mu,D_\omega}\} \\ 
&= \|(a)_{ik}\|_n \|(b)_{pq}\|_n
\end{split}
\]
Hence, by Theorem 5.2.9 of \cite{BleLeM}, $\Ai{n}$ with thus defined matrix norms is $cb$-isomorphic to an operator algebra. We also show explicitly in \cite{Iva0} that this operator algebra may be chosen to be involutive. 

By the construction we also have that $\Ai{n}$ is a subalgebra of $\Ai{n-1}$, and that 
$$\|(a)_{ik}\|_{n,\mu,D} \le \|(a)_{ik}\|_n$$
so that the inclusion map $\Ai{n} \to \Ai{n}_{\mu,D}$ is actually completely contractive.
\end{proof}
\end{lem}

\begin{lem}\label{smth_alg_lem_general_lemma_2}
Let $\mu$ be a commutator bounded analytic fr{\'e}chetization. Then for any separable $C^*$-algebra $A$ and any set $\{(E_\omega,F_\omega)\}_{\omega\in\Om}$ of $\,\KK$-cycles over $(A,B_\omega)$ there exists an $\iy$-smooth system $\Asm$ such that the map $\Qin{\iy}(A,B_\omega)\to\{[(E_\omega,F_\omega)]\}_{\omega\in \Om}$ is surjective.
\begin{proof}
Without loss of generality we may suppose $F^2=1$ and $F=F^*$. Choose a total system $\{a_j\}$ on $A$ and construct the unbounded $\KK$-cycles $(E_\omega, D_\omega)$ for each $(E_\omega,F_\omega)$ by the method described in Lemma \ref{smth_alg_lem_tech_lemma}. Since $\mu$ is commutator bounded, we have that for the algebra $\Ai{n}_{\mu,D}$ holds 
\begin{equation}\label{smth_alg_eqn_aj_extimate}
\|a_j\|_{n,\mu,D} \le C_n\max\{\|a\|, \|\ad_D(a)\|, \dots, \|\ad_D^n(a)\|\} \le C_n \max_{k=0,\dots, n}(c_{k,j}) =: K_{n,j}
\end{equation}
and $K_{n,j}$ does not depend on $\omega$. Thus we may apply the Lemma \ref{smth_alg_lem_general_lemma_1}. We denote the resulting sequence of algebras $\Asm = \{\Ai{n}\}$. To prove that $\Asm$ is a smooth system, we need to prove the holomorphic stability of the algebras $\Ai{n}$. But since $\mu$ is analytic  we have that $\|\cdot\|_{n,\mu,D_\omega}$ is analytic with respect to $\|\cdot\|_{n-1,\mu,D_\omega}$ for all $n$. Observing that $\|a\|_n \le \|a\|_{n,\mu,D_\omega}$ for all $n$, we have that $\|a\|_{n-1,\mu,D_\omega} \le 1$ for all $a$ such that $\|a\|_{n-1}\le 1$, and so for all such $a$ holds

$$
\limsup_{m\to\iy}\frac{\ln\|a^m\|_n}{m} = 
\limsup_{m\to\iy} \frac{\ln(\sup_{\omega\in\Om}\|a^m\|_{n,D_\omega})}{m}=
\limsup_{m\to\iy} \frac{\sup_{\omega\in\Om}\ln\|a^m\|_{n,D_\omega}}{m}\le 0
$$
Thus, the norm on $\Ai{n}$ is analytic with respect to $\Ai{n-1}$, and so $\Ai{n}$ are stable under holomorphic functional calculus. The holomorphic stability of $\Ai{\iy}$ follows immediately from its definition.

To finish the proof we only need to observe that by the construction $[(E_\omega,\gtb(D_\omega))]=[(E,F)]$.
\end{proof}
\end{lem}

The Theorem \ref{smth_alg_thm_relation_to_KK} then becomes an easy corollary of the Lemma \ref{smth_alg_lem_general_lemma_2}. 

\begin{cor}[Proof of Theroem \ref{smth_alg_thm_relation_to_KK}]
\begin{proof}
Indeed, let $\Lambda$ be the set of (isomorphism classes) of $C^*$-algebras. For any $B_\la\in \Lambda$ choose a set $\Omega_\la$, consisting of the $\KK$-cycles $(E_{\lambda_\omega},D_{\lambda_\omega})$, such that the map $\Lambda \to \KK(A,B_\lambda)$, given by taking the homotopy class, is surjective. Then $\Omega := \bigcup_{\la\in\Lambda}\Om_\la$ is a set. Applying the Lemma \ref{smth_alg_lem_general_lemma_2} we obtain the desired result.
\end{proof}
\end{cor}

\begin{rem}
Some conditions of the Theorem \ref{smth_alg_thm_relation_to_KK} may be either relaxed or replaced by other ones. 

For instance, as we have seen from the proof of Lemma \ref{smth_alg_lem_general_lemma_1}, the condition of the analyticity of the fr{\'e}chetization is a technical issue, allowing us to prove that the constructed algebra $\Ai{n}$ is stable under holomorphic functional calculus. Thus, if some other condition on $\mu$ gives the same result, it may freely replace analyticity. As an example one may require the sequence of norms $(\|\cdot\|, \|\cdot\|_{1,\mu,D}, \|\cdot\|_{2,\mu,D}, \dots)$ to be a \emph{differential seminorm} (see, \cite{BlaCu}, \cite{Bha} for definitions), which in particular will imply the holomorphic stability of all $\Ai{n}_{\mu,D}$ in case when $\ord\Asm_{\mu,D} = \iy$. With some additional steps in the construction of algebras $\Ai{n}$ we may then make the sequence $(\|\cdot\|,\|\cdot\|_1, \|\cdot\|_2, \dots)$ to be a differential seminorm as well, proving thus that $\Ai{n}$ are pre-$C^*$-algebras as required.

We may also somewhat relax the condition that the algebras $\Ai{n}_{\mu,D}$ are operator algebras. Namely, we may demand the algebras $\Ai{n}_{\mu,D}$ to be $K_n$-operator algebras in the sense of \cite{BleLeM}, that is, the operator spaces, which are also algebras with the multiplication inducing a completely bounded bilinear map $m\colon \Ai{n}\times \Ai{n} \to \Ai{n}$, such that the $cb$-norm of $m$ is $\le K_n$. 
\end{rem} 

\section{Fr{\'e}chetizations and Their Properties}
As it has been said above, there have already been defined several fr{\'e}chrtizations. We will briefly describe two of them.

For the simplest one, let 
\[
\Te^1_D(a)=
\begin{pmatrix}
a & 0 \\
[D;a] & a
\end{pmatrix}
\]
as an operator on $E^{\oplus 2}$. Analogously we set
\[
\Te^n_D(a)=
\begin{pmatrix}
\Te^{n-1}_D(a) & 0 \\
[D;\Te^{n-1}_D(a)] & \Te^{n-1}_D(a)
\end{pmatrix}
\] 
on $E^{\oplus 2^n}$, where we abusively denote by $D$ the operator $\diag(\underbrace{D,\dots,D}_{2^{n-1}\text{ times}})$ on $E^{\oplus 2^{n-1}}$. We also have that $\Te(ab)=\Te(a)\Te(b)$

Now, define an algebra $\Ai{n}:=\{a\in A\mid \ad^k_{D} \text{ extends to bounded on }E \text{ for all } k=1,\dots, n\}$. This is an algebra, which is complete with respect to the norm $\|\cdot\|_{n,D}:=\|\Te^n_D(a)\|$. By the construction we have an estimate 
$$
\|a\|_{n,D} \le 2^n\max\{\|a\|,\|\ad_D(a)\|,\dots,\|\ad_n(a)\|\}
$$
so that the fr{\'e}chetization is commutator bounded. In case when $\Ai{n}_D$ is dense in $A$, for all $a\in\Ai{n-1}$ such that $\|a\|_{n-1,D}<1$ we have that
\[
\begin{split}
\limsup_{m\to\iy}\frac{\ln \|a^m\|_{n,D}}{m}&=
\limsup_{m\to\iy}\frac{\ln \|\Te^n_D(a^m)\|}{m}\\
&\le \limsup_{m\to\iy}\frac{\ln (\|\Te^{n-1}_D(a^m)\| + \|[D,\Te^{n-1}_D(a^m)]\|)}{m}\\
&\le \limsup_{m\to\iy}\frac{\ln(1+m\|\Te^{n-1}_D(a^m)\|)}{m}\\ 
&\le \limsup_{m\to\iy}\left(\frac{\ln m}{m} +\frac{\ln(1+\|\Te^{n-1}_D(a^m))\|}{m}\right)\\
&=0 
\end{split}
\]

Hence, we are in the conditions of the theorem \ref{smth_alg_thm_relation_to_KK}.

However, there is nothing that may guarantee us that the smooth system $\Asm_D$ will have the order of smoothness exceeding $1$ (order 1 is given by the fact that $(E,D)$ is an unbounded $\KK$-cycle). Moreover, this fr{\'e}chetization is in fact a na{\"i}ve one and may be useful only for theoretical needs lying out of th differential geometry. It is easy to see that even for the Dirac operator on a 2-torus $[D,[D,\cdot]]$ does not extend to a bounded operator on the spinor bundle of the torus. Therefore the fact that $\ad^n_D(a)$ extends to a bounded operator on $E$ is a dense subalgebra of $A$ tells us that the Theorem \ref{smth_alg_thm_relation_to_KK} is entirely an existence result. 

The calculations show the following fact. Take an element $(\Lc^2(S^1),\Dn:=-i\frac{\partial}{\partial x})$ as an unbounded $(C(S^1),\Cb)$-$\KK$-cycle and Fourier functions $\{e^{2\pi i k}\}_{k\in\Zb}$ for a total system as a starting point. Perform the series of operations
$$
\Dn \to \btr{\Dn} \to F \to D 
$$
where the second arrow is a modification by a compact operator, yielding a selfadjoint $F$ with $F^2$=1, and the third is provided by the construction from Lemma \ref{smth_alg_lem_tech_lemma}. Then the sequence $|\la_l|$ of absolute values of eigenvalues of $D$ ordered by increasing will grow at most logarithmically. Moreover, for any $n$ there exists a function $a\in C(S^1)$ which is nowhere differentiable in common sense, but still $\ad^n_D(a)$ extends to a bounded operator on $\Lc_2(S^1)$, so that $D$ regards $a$ as a "$C^n$" function. This fact is shown explicitly in \cite{Iva0}. 

It is also possible to start with a total system $\{a_j\}$ such that all $a_j$ will be nowhere differentiable funcltions, but the result \ref{smth_alg_thm_relation_to_KK} would still be valid.

In \cite{Mesl} there have been constructed a more elaborated fr{\'e}chetization, for which the Dirac operators on spin manifolds do yield $\iy$-smooth systems. It is not called a fr{\'e}chetization there, but is regarded as a natural construction of a smooth system by an unbounded operator. There is also shown in \cite{Mesl} that for this fr{\'e}chetization the operators $D$, $cD$ and $D+b$, with $b$ being a bounded operator satisfying some additional smoothness conditions, generate equivalent smooth systems in the sense that the algebras $\Ai{n}_D$, $\Ai{n}_{cD}$ and $\Ai{n}_{D+R}$ are $cb$-isomorphic. This observation may be crucial for the development of further theory. Indeed, because of that we have that if $(E,D)\in\Qin{n}(\Asm,B)$, then so do $(E,cD)$ and $(E,D+b)$. From the other hand, by definition if a smooth system $\Asm'$ is equivalent to $\Asm$ in the above sense, then we have that $(E,D)\in \Qin{n}(\Asm',B)$. Recall (\cite{ConNCG},\cite{Vari}) that in the case of spectral triples with commutative algebras the Dirac-type operators define metrics on the underlying topological space. Therefore $\Qin{\iy}(\Asm,\Cb)$ in this case may be regarded as the set of metrics that are "smooth" with respect to some system given by $\Asm$. Thus, $\Asm$ in some sense becomes an analogue of the system of smooth functions which in geometry is standardly achieved by imposing coordinates on local charts and taking partial derivations. Thus, the notion of smooth system may lead us to a notion, which may be called a \emph{smooth noncommutative topology}. This question is discussed in more detail in \cite{Iva0}.

\end{document}